\documentclass{amsart}

\usepackage{amssymb}
\usepackage{amsmath}

\theoremstyle{definition}

\theoremstyle{remark}

\numberwithin{equation}{section}

\begin{document}


\title{Some Technical Thoughts on Modeling}


\author{Nikolas O. Aksamit$^1$}
\address{$^1$Department of Geography and Planning, University of Saskatchewan, Saskatoon, SK S7N 5C8}
\email{n.aksamit@usask.ca}



\author{Don H. Tucker$^2$}
\address{$^2$Department of Mathematics, University of Utah, Salt Lake City, UT 84112}

\author{James F. Tucker$^3$}
\address{$^3$School of Medicine, University of Utah, Salt Lake City, UT 84112}




\begin{abstract}

This paper is more an essay than a report. There is a gentle introduction to some issues in modeling, followed by the use of steepest descent methods to develop a model as contrasted to using such methods to solve one already in hand, as in \cite{KA}. Three levels are discussed: fitting functions to model given data, fitting an ODE to model given data, and more briefly, fitting a PDE to model given data. Specific examples are discussed.

\end{abstract}


 \maketitle




AMS Keywords: 00-02 (Research Exposition), 92-08 (Computational Methods), 65C20 (Models, Numerical Methods)

\section{Introduction}

Concepts are not born fully mature on the half shell like a Venus as in Botticelli's rendition of Zeus' spit. Concepts develop slowly from small fuzzy shadows of awareness. They do not reveal themselves to any but those who are alert to their presence and who work hard to place themselves in positions from which the concepts can be seen. Such work nearly always involves personal experience which is tediously examined from various vantage points of thought. Such experience is frequently the collective personal encounters of the species over extended periods of time together with the recorded thought examinations communicated from one generation to the next. Some have suggested that the major importance of the invention of writing has been the warehousing of our thoughts for the use of those who come later. Our efforts to model the universe are of this collective nature.

The bulk of our efforts at modeling are devoted to local models; models whose scope is restricted to a limited set of events in a small amount of space for a short period of time. We then extrapolate to other cases, more extensive in space and time.

Some of the ways we do such things include the following. We observe some event, measure some aspects of it which we can measure, conjecture some parameters which might influence the things happening and then conjecture functions of those parameters which (hopefully) will result in the measurements we have made. We assume that these parameters are "physical" parameters in the sense that they are not time dependent; if our model is valid  today, it will be valid next week as well.

Assume $x(t)$ is one of the attributes we measure over our time interval at discrete times $\{x(t_i)\}_{i=1}^{n+1}$, $a=(a_1, ..., a_p)$ are our conjectured physical parameters, and $f(a,t)$ is our model of $x(t)$. The conjectured parameters $a$ are not what we measured. We have a strong desire to have rather precise estimates of those parameters which we may then insert into the function $f$ and if this replicates $x(t)$ with adequate precision, we extrapolate beyond our observations.

In the event that our only data is the set of our measurements $\{x(t_i)\}$ and we are (for whatever reasons) limited to those data, we usually resort to steepest descent methods to obtain estimates for the parameters $a$. We address this in section 2.

The basic mathematical idea we use is not new. It was discussed by Cauchy for finite dimensional spaces and extended to more general spaces by Kantorovich and others as in \cite{KA}. The use was to find solutions to a given problem by reducing the problem to a variational problem and using minimizing techniques to achieve a result. Our intention here is to start with the solution, that is observed data, and find a problem, that is a model, whose solution will be the observed solutions. We are not so naive as to think the result, if achieved, will always be unique.

In other situations our conjectures may involve rates of change of the measured data. The Tyco Brahe, Kepler, Galileo case of planetary observations, conjectures and calculations were of this nature and likely led to the the result of Isaac Barrow (the first part of the fundamental theorem of calculus) which he taught to Isaac Newton. The situation is roughly like this: we have measurements $\{x(t_i)\}$, we make conjectures, $x'(t)=f(a, x(t))$, and require (desire) the parameters $a$ in the resulting differential equations. This is the concern of section 3.

If we believe our observations are the result of rates of change in more than one independent dimension, our model will involve partial derivatives in multiple dimensions. We address the method of using given data to approximate the governing PDE's of our system in section 5.

\section{Watersheds}

Fitting a conjectured function $f(a, x)$ to a given set of data $\{x(t_i)\}$ involves making a decision as to what measure one will use to define "fit." In this note we will mean that the sums of the squares of the discrepancies $|f(a, x(t_i))-x(t_i)|$ 
is made small enough to comport with the prescribed precision. Thus, our initial goal is to find an $a$ such that

$$F(a)=\sum_{i=1}^n|f(a, x(t_i))-x(t_i)|^2<m,$$

then certainly $|f(a, x(t_i))-x(t_i)|<\sqrt{m}$ at each datum $x(t_i)$, where $\sqrt{m}$ is the desired precision at each of the individual points.

The method used is: Guess a starting point, $a_1$. Compute $F(a_1)$ and the gradient of $F(a)$ at $a=a_1$. The maximum rate of decrease of $F$ as a function of $a$ is in the direction of $-\nabla F|_a$. Guess a value $\epsilon$ which will be our step size from $a_1$ to $a_2=a_1-\nabla F|_{a_1}\cdot\epsilon$. If $F(a_2)$ is smaller than $F(a_1)$, continue. That is, start at $a_2$ and have $a_3=a_2-\nabla F|_{a_2}\cdot\epsilon$. If ever $F(a_{n+1})>F(a_n)$, back up and try a smaller $\epsilon$.

In as much as the function $F(a)$ may have more than one minimum point, different choices of $a_1$ could result in multiple sets of parameters $a$ which meet one's requirements. If  $m$ is a point where $F$ attains a local minimum, we define the watershed of $m$ as the neighborhood of $m$ inside which continuous flow in the direction of the negative gradient will lead to $m$. If the initial guess $a_1$ is within the watershed of a local minimum, and a single step of size $\epsilon$ does not escape the watershed, then a least squares regression can be trapped, unable to escape to find other more global minima.  To approach this problem, one can use what we call the "shotgun" approach, where many initial guesses are made at a variety of locations in the parameter space, the minimization results of which can then be compared after the fact to select the best candidate.  Other, more sophisticated methods, rely on using different values of $\epsilon$ at different times in the process in an attempt to "jump" out of such watershed traps.  See \cite{Kraft}.

The second difficulty with this method is determining the reliability of the fit in the face of experimental noise in $\{x(t_i)\}$.  This problem is not unique to the least squares method, nor is it specific even to the fitting of functions.  All models derived from real-world data must be carefully examined for the extent to which major features of the generated model are sensitive to small changes in the initial data, lest the model fail to describe the general case. For an example of generating models from real-world data using functional fitting, including noise stability testing, see \cite{Tucker}.

Those who have encountered statistics will recognize this procedure as linear regression in the event that $f(a, x)$ is assumed to be a linear map in $x$. The procedure has been used to considerable advantage in \cite{Tucker} in cases in which the raw data is thought to have been generated by multiple simultaneous processes which are each represented by Gaussian distributions. Separating such distributions then led to better understanding of the phenomena involved.

\section{ODE Estimate}

Suppose we have a set of points $(t_i, x(t_i))$ in $\mathbb{R}^{n+1}$ space, i.e., $x(t_i)\in\mathbb{R}^n$ and we wish to construct an ODE, $x'(t)=f(t,x(t))$ whose solutions (which satisfy given initial data) replicate the above data points to within some precision yet to be determined. How close can we come, whatever that means?

Our method (of madness) is as follows: We guess a function $f(a, t, x(t))$ which is reasonably smooth (we will assume $C^{(1)}$ as we proceed) and which might come close if the parameters $a$ are suitably chosen. That is, we conjecture a model of the situation being observed. We will use steepest descent methods to determine an acceptable set of $a$'s once the function $f$ has been conjectured. Guessing the $f$ will almost certainly (not a probabilistic term) depend upon the past experience of the guesser with the phenom which produced the data points $\{t_i, x(t_i)\}$. We have no advice nor algorithms to offer in this regard. Bridgman may have said it best, "The problem cannot be solved by the philosopher in his armchair, but the knowledge involved was gathered only by someone at some time soiling his hands with direct contact." \cite[p. 11-12]{Bridgman}

Set $F(a)=\sum_{i=1}^n \|f(a,t_i, x(t_i))-\frac{[x(t_{i+1})-x(t_i)]}{t_{i+1}-t_{i}}\|^2$ and minimize $F(a)$ by steepest descent as a function of $a$. Assume $f$ is $C^{(1)}$ in the parameters $a$.

Suppose that is done, $a$ is determined and thus $f$ is fixed so that

$$\sum_{i=1}^n \|f(t_i, x(t_i))-\frac{[x(t_{i+1})-x(t_i)]}{t_{i+1}-t_{i}}\|^2=\sum_{i=1}^n \|x(t_i)+f(t_i, x(t_i))[t_{i+1}-t_i]-x(t_{i+1})\|^2\frac{1}{(t_{i+1}-t_i)^2}\le m,$$ where $m$ is the minimum value achieved by steepest descent.

Define $$p(t)=x(t_i)+f(t_i, x(t_i))[t-t_i]\text{; } t_i\le t< t_{i+1}\text{ with } p(t_1)=x(t_1)$$ and $$P(t)=x(t_i)+\frac{[x(t_{i+1})-x(t_i)]}{t_{i+1}-t_i}[t-t_i]\text{; } t_i\le t\le t_{i+1}.$$ $P(t)$ is a polygonal function whose graph connects the successive data points and thus $P(t)$ is continuous. However, for $p(t)$ we have $$\lim_{t\to t_i^{+}}p(t)=x(t_i) \text{ but} \lim_{t\to t_i^{-}}p(t)\neq x(t_i).$$ Let us call $\Delta t_i=|t_{i+1}-t_i|$. Rewriting what we have above,

$$\sum_{i=1}^n\|p(t_{i+1})-x(t_{i+1})\|^2\frac{1}{(\Delta t_i)^2}\le m.$$ Suppose $A\le \Delta t_i \le B$ for $i=1,2, ..., n$, then

$$\sum_{i=1}^n\|p(t_{i+1})-x(t_{i+1})\|^2\le \sum_{i=1}^n\|p(t_{i+1})-x(t_{i+1})\|^2\frac{B^2}{(\Delta t_i)^2}\le mB^2$$ and therefore $\|p(t_{i+1})-x(t_{i+1})\|\le \sqrt{m}B$ for each $i$. This is a large overestimate, but the best we can afford currently. It is worth noting

$$\|p(t)-P(t)\|=\|f(t_i, x(t_i))-\frac{[x(t_{i+1})-x(t_i)]}{[t_{i+1}-t_i]}[t-t_i]\|\le\sqrt{m}B\Delta t_i\le \sqrt{m}B^2$$ and this is uniform over the entire $t$ domain.

Recall that $p'(t)=f(t_i, p(t_i)); t_i\le t< t_{i+1}$ and $p(t_1)=x(t_1)$. With the $f$ now determined via steepest descent, consider a solution $y'(t)=f(t, y(t)); y(t_1)=x(t_1)$.

Our concern (just now) is how small is $\|P(t)-y(t)\|$ on the domain of $t$. The pursuit of an answer is by way of $p(t)$ since we already have a measure $\|p(t)-P(t)\|\le\sqrt{m}B^2$ for every $t$.

\begin{align*}
p'(t)-y'(t)&=f(t_i, p(t_i))-f(t,y(t))
\\&=f(t,p(t))-f(t,y(t))+f(t_i,p(t_i))-f(t,p(t))
\end{align*}
$\text{for } t_i<t\le t_{i+1}$.

Set $g(t)=f(t_i, p(t_i))-f(t,p(t))$ for $t_i\le t< t_{i+1}$ and integrate from $t_1$ to $t$. There is no loss in assuming $t_1=0$; let's do that. (Remark: We know almost nothing about $f(t, p(t))$ without further assumptions).

$$[p(t)-y(t)]-[p(t_1)-y(0)]\le\int_{0}^t\|g(u)\|du+\int_{0}^t\|f(u, p(u))-f(u,y(u))\|du$$ Assuming $f$ is Lipschitz on its domain with some Lipschitz constant $L$, we get

$$\|f(u,p(u))-f(u,y(u))\|\le L\|p(u)-y(u)\| \text{ hence we have}$$
$$\|p(t)-y(t)\|\le\int_{0}^t\|g(u)\|du+L\int_{0}^t\|p(u)-y(u)\|du$$ since $p(t_1)=x(t_1)=y(t_1)$. Now, set $F_m=\text{max}_i\{\|f(t_i, x(t_i))\|\}$. At this point, we wish to find an upper bound for $\|g(t)\|$, namely

\begin{align*}
\|f(t_i, p(t_i))-f(t,p(t))\|&\le L\|(t,p(t))-(t_i, p(t_i))\| 
\\ &\le L[(t-t_i)^2+(p(t)-p(t_i))^2]^{\frac{1}{2}}
\\ &\le L[(t-t_i)^2+(t-t_i)^2F_m^2]^{\frac{1}{2}} \le LB[1+F_m^2]^{\frac{1}{2}}
\end{align*}

We'll give this upper bound for $\|g\|$ a name, say $M$. Then we have that

\[\|p(t)-y(t)\|\le Mt+L\int_{0}^t\|p(u)-y(u)\|du.\] We need a modified form of Gronwall's inequality. Suppose $f\ge0$ and $g\ge0$ on $0\le t$ and $M>0$ and $f(t)\le Mt+\int_0^tf(s)g(s)ds$.

Set $$H(t)=Mt+\int_0^tf(s)g(s)ds,$$ then $f(t)\le H(t)$ and $$H'(t)=M+f(t)g(t)\le M+[Mt+\int_0^tf(s)g(s)ds]g(t)\le M+H(t)g(t)$$ Multiply by $e^{-\int_0^tg(s)ds}$, an integrating factor, and get

\[
H'(t)e^{-\int_0^tg(s)ds}\le[M+H(t)g(t)]e^{-\int_0^tg(s)ds}\tag{$\dagger$}
\]

Notice that 

\begin{align*}
[H(t)e^{-\int_0^tg(s)ds}]'&=H'(t)e^{-\int_0^tg(s)ds} + H(t)[-g(t)]e^{-\int_0^tg(s)ds}
\\&=H'(t)e^{-\int_0^tg(s)ds}-H(t)g(t)e^{-\int_0^tg(s)ds}
\end{align*}

Hence, $(\dagger)$ becomes

$$H'(t)e^{-\int_0^tg(s)ds}-H(t)g(t)e^{-\int_0^tg(s)ds}\le Me^{-\int_0^tg(s)ds}$$

\[
\text{or   } [H(t)e^{-\int_0^tg(s)ds}]'\le Me^{-\int_0^tg(s)ds}\tag{$\ddagger$}
\]

Integrate both sides from $0$ to $t$ and get

$$H(t)e^{-\int_0^tg(s)ds}-H(0)\le M\int_0^te^{-\int_0^ug(s)ds}du$$ or

$$H(t)\le\{M\int_0^te^{-\int_0^ug(s)ds}du\}e^{\int_0^tg(s)ds}$$

\hfill$\Box$

Upshot: In our case $g(u)\equiv L$, the Lipschitz constant, which gives:

\begin{align*}
f(t)&\le\{M\int_0^te^{-Lu}du\}e^{Lt}=Me^{Lt}\int_0^te^{-Lu}du
\\&=Me^{Lt}\frac{-1}{L}[e^{-Lt-1}]=\frac{M}{L}e^{Lt}[1-e^{-Lt}]=\frac{M}{L}[e^{Lt}-1]
\end{align*}
and $f(t)=\|p(t)-y(t)\|$ and $M=\|g\|$. Recall $M\le LB[1+F^2]^{\frac{1}{2}}$. This gives the result:

\begin{align*}
\|p(t)-y(t)\|\le\frac{M}{L}[e^{Lt}-1]&\le\frac{LB[1+F^2]^{\frac{1}{2}}}{L}[e^{Lt}-1]
\\ &= B[1+F^2]^\frac{1}{2}[e^{Lt}-1]
\end{align*} for every $t\ge0$

It follows as night the day that

\begin{align*}
\|P(t)-y(t)\|&\le\|P(t)-p(t)\|+\|p(t)-y(t)\|
\\
&\le \sqrt{m}B^2+B[1+F^2]^{\frac{1}{2}}[e^{Lt}-1]
\\
&= B\{\sqrt{m}B+[1+F^2]^{\frac{1}{2}}[e^{Lt}-1]\}
\end{align*}

This is small provided $B$, $m$ and $[e^{Lt}-1]$ are small. The first two require precision of measurements and calculations while the third requires that $t$ be near zero. This shows that our model is local in nature from the mathematical structures involved, not just from the physical considerations mentioned above. This also indicates that perturbations in precision possibly propagate quite rapidly. If $t$ is measuring time or if $t$ is measuring distance, one is cautioned just the same; reliability may well degrade as the model is pushed farther. 

Question: Can we estimate $F=max\|f(t_i, x(t_i)\|$ directly from $\{x(t_i)\}$'s and $m$? If so, our error estimates would be almost independent of the choice of $f$, but would depend on $L$ and $m$. We shall say $max\|\frac{\Delta x_i}{\Delta t_i}\|=\Delta$. Notice:

\begin{align*}
\sum\|f(t_i,x(t_i))-\frac{\Delta x_i}{\Delta t_i}\|^2 & \le m
\\
\Rightarrow \|f(t_i,x(t_i))-\frac{\Delta x_i}{\Delta t_i}\| & \le\sqrt{m}
\end{align*}
$$\|f(t_i,x(t_i))\|\le\|\frac{\Delta x_i}{\Delta t_i}\|+\sqrt{m}$$

$$F\le\sqrt{m}+max\|\frac{\Delta x_i}{\Delta t_i}\|=\sqrt{m}+\Delta$$

$$\|P(t)-y(t)\|\le B\{\sqrt{m}B+[1+(m+\Delta)^2]^{\frac{1}{2}}[e^{Lt}-1]\}$$

If we restrict $L$ (physically this is restricting $y''$, the acceleration or force) and require $m$ be smaller than a certain fixed precision, we may be able to give a comparison result between solutions which result from different models, each derived from the same data by these methods.



Suppose Joe Blow conjectures a different $C^{(1)}$ function $h$, rather than $f$. The steepest descent methods afford him a total error of $\tilde{m}$. The maximum for $h$ over our domain is $H$, and $h$ has a Lipschitz constant $\tilde{L}$. Note that $P(t)$ is the same for Joe as it is for us.

Joe then gets a solution $z(t)$ to his ODE, $z'(t)=h(t,z(t))$, $z(0)=x(t_1)$. How different are our conjectured models of reality?

$$\|P(t)-z(t)\|\le B[\sqrt{\tilde{m}}B+\tilde{L}(t-t_1)[1+H^2]^\frac{1}{2}e^{\tilde{L}(t-t_1)}]$$

$$\|P(t)-y(t)\|\le B[\sqrt{m}B+L(t-t_1)[1+F^2]^\frac{1}{2}e^{L(t-t_1)}]$$

$$\therefore \|y(t)-z(t)\|\le B\{(\sqrt{m}+\sqrt{\tilde{m}})B+[1+F^2]^{\frac{1}{2}}[e^{Lt}-1]+[1+H^2]^{\frac{1}{2}}[e^{\tilde{L}t}-1]\}$$

If $m$ and $\tilde{m}$ are required to be $<\delta^2$, we have $F$, $H<\delta+\Delta$ in which case

$$\|y(t)-z(t)\|\le B\{2\delta B+[1+(\delta+\Delta)^2]^{\frac{1}{2}}[(e^{Lt}-1)+(e^{\tilde{L}t}-1)]\}.$$

If $L$, and $\tilde{L}$ are required to be $<\mathcal{L}$, we would have

\begin{align*}
\|y(t)-z(t)\|&\le B\{2\delta B+[1+(\delta+\Delta)^2]^{\frac{1}{2}}[2(e^{\mathcal{L}t}-1)]\}
\\
&\le 2B[\delta B+[1+(\delta+\Delta)^2]^{\frac{1}{2}}[e^{\mathcal{L}t}-1]
\end{align*}
for all "acceptable" solutions $y$ and $z$, where $B$ and $\Delta$ are determined by the raw data and $\delta$ and $\mathcal{L}$ are imposed for physical reasons. The right hand side then is independent of the choices of $f$ and $h$. This brings us to a vexing but quite real scientific problem. Suppose additional data points are not to be had, for whatever reasons. Further assume there are several models which give acceptable precision at the data points yet differ greatly if extended much beyond the initial local domains for $t$ and $x$. A major reason for building a model is to use it to predict beyond the observed situation. If the different models predict differently, how does one choose among them short of more observations? Again, we have no advice concerning a choice among such models.

Before we compute an example, let us consider the experimental noise mentioned earlier. Suppose there is a potential measurement error $\epsilon$ for each data point $x(t_i)$ such that there is a flag; $x(t_i)\pm \epsilon_i$ and $\epsilon=max |\epsilon_i|$. Then we bound each $x(t_i)$

$$\underline{x}(t_i)=x(t_i)-\epsilon\le x(t_i) \le x(t_i) + \epsilon =\bar{x}(t_i).$$

Repeating our steepest descent method for $\{\underline{x}(t_i)\}$ and $\bar{x}(t_i)$ we obtain functions $\underline{f}$ and $\bar{f}$, respectively, as well as an $f$ for our measured data points, $\{x(t_i)\}$. Suppose these functions then give rise to solutions, $\underline{y}$, $\bar{y}$ and $y$. Assume we solve for $f$ first, and use the resulting parameter point $a$ as the initial guess to solve for $\underline{a}$ and $\bar{a}$. Hopefully these remain in the same watershed. Also, mutatis mutandis, denote $\underline{m}$, $\bar{m}$, $m$, $\underline{P}$, $\bar{P}$, $P$, $\underline{L}$, $\bar{L}$, $L$, and $\underline{F}$, $\bar{F}$, $F$. Let's compare $\bar{y}$ and $y$ (change notation to compare $\underline{y}$ and $y$).

First note that $B$ is the same for all cases, as it is irrespective of measurement error, and $|\bar{P}(t)-P(t)|\le\epsilon$. We now have that

\begin{align*}
\|\bar{y}(t)-y(t)\|&\le\|\bar{y}(t)-\bar{P}(t)+\bar{P}(t)-P(t)+P-y(t)\|
\\
&\le\bar{y}(t)-P(t)\|+\epsilon+\|P(t)-y(t)\|.
\end{align*}

We already have bounds for $\|\bar{y}(t)-\bar{P}(t)\|$ and $\|P(t)-y(t)\|$, which admittedly may be gross overestimates. Nonetheless, these give us the impact of our measurement error:

\[
\|\bar{y}(t)-y(t)\|\le\epsilon+\|\bar{y}(t)-P(t)\|+\|P(t)-y(t)\|=\epsilon+\bar{E}+E
\]
\\
where 

\begin{align*}
\bar{E}&\le B[\sqrt{\bar{m}}B+\bar{L}(t-t_1)[1+\bar{F}^2]^{\frac{1}{2}}]e^{\bar{L}(t-t_1)} \text{ and}
\\
E&\le B[\sqrt{m}B+L(t-t_1)[1+F^2]^{\frac{1}{2}}]e^{L(t-t_1)}
\end{align*}

This is similarly done with $\underline{y}(t)$ to then obtain the size of the entire neighborhood of error: 

\[
\|\bar{y}(t)-\underline{y}(t)\|=\|\bar{y}(t)-y(t)+y(t)-\underline{y}(t)\| \le 2(\epsilon+E)+\bar{E}+\underline{E}
\]




\section{Computed Examples}

Heeding Bridgman's remarks concerning soiling one's hands, we checked our methods against several ODEs. The computations were done using Matlab. We first sought to replicate the coefficients in the equation

$$x'(t)=x(t)^2+2x(t)$$ with the initial condition x(0)=1. It has solution $x(t)=\frac{-2e^{2t}}{e^{2t}-3}$. Using values of that known solution for the $x(t_i)$ data and assuming $f(a,x)=a_1x^2+a_2x$, we minimized $$F(a)=\sum_{i=1}^n|a_1x(t_i)^2+a_2x(t_i)-\frac{x(t_{i+1})-x(t_i)}{t_{i+1}-t_i}|^2$$ using steepest descent on several domains for the $\{t_i\}$.

Notice $x(t)$ has a singularity at $\frac{ln(3)}{2}\approx.549$. Away from this point, for example, $1\le t_i\le 2$, with a uniform $\Delta t_i$ as coarse as $\frac{1}{10}$ we were able to retrieve $[a_1, a_2]=[1.00, 2.00]$ with a gradient of $F=[10^{-11}\times.3638, 0]$ and value of $F$ as small as $3.589\times10^{-20}$. Using these values of $a_1, a_2$ we, of course, exactly replicated our "observed" data. 

When we included $\frac{ln(3)}{2}$ as an interior or boundary point in the domain of $t$, the desired Lipschitz condition on $f$ was no longer satisfied because $x(t)$ was unbounded and our errors suffered. For example, with $\Delta t_i=\frac{1}{1000}$, $t_1=0$ and $t_{999}=1$, an initial guess of $[a_1, a_2]=[1.00,2.00]$, gave a value of $F =2.3596\times10^{14}$. This is rather startling considering that we started with exactly the correct values for $a_1$ and $a_2$. In this case, our steepest descent method converged to $[a_1, a_2]=[-.0003\times10^3, -1.9972\times10^3]$ with the value of $F$ approximately $1.4288\times10^{13}$. The values for the actual solutions $y$ of the resulting differential equation differed from the original data, $x$, as follows: 
$$(\sum_{i=0}^{999}|y(t_i)-x(t_i)|^2)^{\frac{1}{2}}\approx3.7341\times10^8.$$ 


Local minima are something that should always be taken into consideration when performing steepest descent, but do not necessarily mean absolute failure of our method. If in our steepest descent computations we use an initial vector guess of $[4,5]=[a_1, a_2]$, and a fixed $\Delta t_i=\frac{1}{1000}$ with this example, then we fall into the local minimum of $[a_1, a_2]=[2.8, 5.6]$. With this information, the norm of the difference in the solution values of the ODE $x'(t)=2.8x^2+5.6x, x(0)=1$, and our initial data is as large as 17.9619 if $t_0=1$ and $t_{999}=2$. For other domains of t, for example, $t_0=10$ and $t_{999}=11$ the norm of our difference is $1.9384\times10^{-7}$, and $t_0=19, t_{999}=20$ gives a difference of $3.4968\times10^{-15}$. We thus have an example of a $t$ domain, an $f=x^2+2x$ and an $h=2.8x^2+5.6x$ where the solutions are indeed quite close on the domain of $t$. This is an example of the issue addressed at the end of section 3. This also illustrates the gross nature of our upper bounds.

\section{PDE Case}

The complicated nature of PDE theory has thus far prevented the comprehensive error inequalities that were possible with ODE's. However, an analog of the ODE steepest descent technique has proven effective at predicting coefficients with several constant coefficient PDE's. First, some notation:

We say the vector $\alpha=(\alpha_1, \alpha_2, \cdots, \alpha_n)$ is of order $|\alpha|=\alpha_1 + \cdots +\alpha_n$. Given a vector $\alpha$ and a differentiable function $u: \mathbb{R}^n\to\mathbb{R}$, we define the differential operator $$D^{\alpha}u(x)=\frac{\partial^{|\alpha|}u(x)}{\partial x_1^{\alpha_1}\cdots\partial x_n^{\alpha_n}}.$$ If after soiling our hands we conjecture our observation data is representative of a function $u$ that satisfies some PDE of the form

$$f(D^mu(x), D^{m-1}u(x), \cdots, Du(x), u(x), x)=\sum_{|\alpha|=0}^m a_{\alpha}(x) D^\alpha u(x)+c=0$$ with specific boundary conditions, where, $c\in\mathbb{R}$ and $a_{\alpha}(x):\mathbb{R}^n\to\mathbb{R}$, then our next task is to solve for the $a_\alpha$. In all of our computed examples, we have worked with the $a_\alpha$ being constant coefficients, but there is no reason to believe our methods would not work with non-linear PDE's. 

After a choice of $f$, we replaced partial derivatives with linear approximations involving our observation data. A general notation for this process would be overwhelming. We will illustrate with some examples. For a function of two variables, a partial derivative in one variable at the point $(x_{i+1}, t_{j+1})$ can be approximated by:

$$\frac{\partial u(x_{i+1},t_{i+1})}{\partial x}\approx \frac{u(x_{i+1}, t_{j+1})-u(x_{i}, t_{j+1})}{x_{i+1}-x_{i}}$$ A second partial derivative at $(x_{i+1}, t_{j+1})$ can be approximated by $$\frac{\partial^2 u(x_{i+1},t_{i+1})}{\partial x^2}\approx \frac{u(x_{i+1}, t_j)-2u(x_{i}, t_j)+u(x_{i-1}, t_j)}{(x_{i+1}-x_{i})(x_i-x_{i-1})}$$ This same process can be extended for mixed and higher order derivatives. It is worth noting that the amount of data you have limits the order of the derivative you can approximate: with observations at $n$ different values in the $x_i$ dimension, you cannot approximate an $n^{\text{th}}$ partial derivative in $x_i$.

To find the coefficients in our PDE $f=0$, we build a new function $F(a, u, x)$ similar to that in our ODE case and minimize in the $a$ dimension via steepest descent.

For one of our examples we used $$u(x,t)=\frac{1}{\sqrt{4\pi 7t}}e^{\frac{-x^2}{4(7t)}},$$ a fundamental solution to the diffusion equation $u_t=7u_{xx}$, as our observed data. We conjectured the PDE was of the form $a_1u_x+a_2u_{xx}+a_3u_t+a_4u_{tt}=0$. Using the range $2\le x_i\le3$, $2\le t_i\le3$, and uniform grid $\Delta x_i=\Delta t_i=\frac{1}{40}$, we minimized

\begin{align*}
F(a,u,x,t)=\sum_{j=1}^n\sum_{i=1}^n &\Big(a_1[u(x_{i+1}, t_{j+1})-u(x_i, t_{j+1})]\frac{1}{\Delta x_{i+1}} 
\\&+a_2[u(x_{i+1}, t_{j+1})-2u(x_i, t_{j+1})+u(x_{i-1}, t_{j+1})]\frac{1}{\Delta x_{i+1}^2}
\\&+a_3[u(x_{i+1}, t_{j+1})-u(x_{i+1}, t_{j})]\frac{1}{\Delta t_{i+1}}
\\&+a_4[u(x_{i+1}, t_{j+1})-2u(x_{i+1}, t_{j})+u(x_{i+1}, t_{j-1})]\frac{1}{\Delta t_{i+1}^2}\Big)^2
\end{align*} using steepest descent. During one computation, we started at the point $[1, -1, 1, 1]$, we were able to obtain $[a_1, a_2, a_3, a_4]=[-.0002, -1.1241, .1631, .0034]$ with the value of $F$ being $1.874\times10^{-10}$. The computations were stopped short of convergence but $a_1$ and $a_4$ were nearing zero, and the ratio $\frac{a_2}{a_3}$ was appearing to be converging to $-7$.

\end{document}